\documentclass[11pt]{article}
\usepackage{latexsym}
\usepackage{amssymb}
\font\gorditas = msbm8
\def\bbb#1{\hbox {{\gordas #1}}}
\def\errita{\hbox{\gorditas R}}
                                                                                
\font\gordas = msbm10 at 12pt
\def\bbb#1{\hbox {{\gordas #1}}}
\def\erre{{\bbb R}}

\def\UNO{1\mkern-7mu1}

\def\flecha{\hbox{$\Rightarrow_{_{_{\!\!\!\!\!\!\!\!\!C}}}\;$}}
\def\flechita{\hbox{$\,\,\Rightarrow_{_{_{\!\!\!\!\!\!\!\!\!f\,\,\,}}}\;$}}

\newtheorem{theorem}{Theorem}[section]
\newtheorem{lemma}[theorem]{Lemma}

\newtheorem{proposition}[theorem]{Proposition}
\newtheorem{ex}[theorem]{Example}

\newtheorem{remark}[theorem]{Remark}

\begin{document}
\vglue.5cm

\begin{center}
{\bf\Large  
Occupation time limits of inhomogeneous\\[.3cm] 
Poisson systems  of  independent particles}
\end{center}
\vglue.3cm
\centerline{T. Bojdecki$^{a,1}$, L.G. Gorostiza$^{b,*,2}$ and 
A. Talarczyk$^{a,1}$ }
\vglue.3cm
\noindent
$^a$ Institute of Mathematics, University of Warsaw, ul. Banacha 2, 02-097 Warsaw, Poland\\
$^b$ Department of Mathematics, Centro de Investigaci\'on y de Estudios Avanzados, A.P. 14-740,\\
\hglue.3cm M\'exico 07000 D.F., Mexico
\footnote{\kern-.6cm$^*$ Corresponding author\\
$^1$ Research partially supported by MNiSW grant 1P03A01129 (Poland).\\
$^2$ Research partially supported by CONACyT grant 45684-F (Mexico).\\
{\it E-mail addresses:} tobojd@mimuw.edu.pl (T. Bojdecki), lgorosti@math.cinvestav.mx (L. G. Gorostiza),\\
 annatal@mimuw.edu.pl (A. Talarczyk).}
\vglue.5cm
\noindent
{\bf Abstract}
\vglue.5cm
We prove  functional limits theorems for the occupation time process of a system of particles moving independently in $\erre^d$ according to a symmetric $\alpha$-stable L\'evy process, and starting off from an inhomogeneous Poisson point measure with intensity measure $\mu(dx)=(1+|x|^{\gamma})^{-1}dx,\gamma>0$, and other related measures. In contrast to the homogeneous case $(\gamma=0)$, the system is not in equilibrium and ultimately it vanishes, and there are more different types of occupation time limit processes depending on arrangements of the parameters $\gamma, d$ and $\alpha$. The case $\gamma<d<\alpha$ leads to an extension of fractional Brownian motion.
\vglue.5cm
\noindent
\vglue.5cm
\noindent
\vglue1cm
\noindent
{\bf 1 Introduction}
\setcounter{section}{1}
\vglue.5cm
Several authors have studied systems of particles moving independently in $\erre^d$ according to a Markov process (usually a symmetric $\alpha$-stable L\'evy process, $0<\alpha\leq 2$), and also systems having in addition a branching mechanism (e.g. \cite{BGT1, BGT2, BGT3, BGT4, CG, DGW, DR, DW, GW, HS, K, ML, Mi, T, W} and references therein). A typical assumption in the cited references  is that the system starts off from a homogeneous Poisson point measure, i.e., with intensity  the Lebesgue measure (denoted here by $\lambda$). This assumption represents a strong technical simplification because in the special cases usually studied $\lambda$ is invariant for the semigroup of the motion, and this implies that the particle system without branching is in equilibrium, and for $d>\alpha$ 
a critical  branching system converges towards equilibrium \cite{GW}. In this case the  systems have been extensively studied. New situations appear if the initial condition is an inhomogenous Poisson point measure.

In this paper we consider the  system without branching, with symmetric $\alpha$-stable L\'evy process for the particle motion, and initial inhomogeneous Poisson point measure with intensity measure $\mu$ of the form
$$\mu(dx)=\frac{dx}{1+|x|^\gamma},\quad \gamma>0,$$
and other more general related measures. In this case the system is not in equilibrium and ultimately it vanishes (see below). Therefore one should expect different types of results from those of the homogeneous case. Our purpose is to obtain  functional limits for the rescaled occupation time process of the particle system in different cases.

The particle system is described as follows.
Given a Poisson point measure on $\erre^d$ with intensity measure $\mu$, particles evolve from its atoms,  moving independently according to a symmetric $\alpha$-stable L\'evy process (called standard $\alpha$-stable process). Let $N=(N_t)_{t\geq 0}$ denote the empirical measure process of the system, i.e.,
\begin{equation}
\label{eq:1.1}
N_t=\sum_i\delta_{x_i(t)},
\end{equation}
where $\{x_i(t)\}_i$  are the positions of the particles at time $t$.
Note that $N_t$ converges in probability to the null measure as $t\to\infty$ 
(Appendix).
 Let $X_T=(X_T(t))_{t\geq 0}$ denote the normalized occupation time fluctuation process of the system, defined by
\begin{equation}
\label{eq:1.2}
X_T(t)=\frac{1}{F_T}\int^{Tt}_0(N_s-EN_s)ds,
\end{equation}
where $T$ is the time scaling and $F_T$ is a norming. The problem  is to find $F_T$ such that the process $X_T$ converges in distribution 
 as $T\to\infty$ (i.e., the time is accelerated), and to describe the limit process $X$ in the cases where it exists. 

In the homogeneous case (corresponding to $\gamma=0$), the occupation time fluctuation limit process has three different forms, for $d<\alpha$ \cite{BGT1}, $d=\alpha$ and $d>\alpha$ \cite{BGT2}. In the inhomogeneous case there are more results depending on the values of $\gamma$ relative to $d$ and 
$\alpha$ when $\gamma\leq d: \gamma< d<\alpha,\gamma<d=\alpha, \gamma<\alpha<d, \gamma=\alpha<d, \gamma=d<\alpha, \gamma=d=\alpha, \alpha<\gamma\leq d$. 
For ``small'' $\gamma$, i.e., $\gamma<\alpha$, the results are analogous to those of  the homogeneous case, while for ``large'' $\gamma$, i.e., $\gamma\geq \alpha$, and this seems  unexpected, they  are of a different kind. 
The case $\gamma<d<\alpha$  leads to a long range dependence, self-similar, centered Gaussian process $\xi$ with covariance
$$E\xi_t\xi_s=\int^{s\wedge t}_0u^{a}((t-u)^b+(s-u)^b)du,$$
where $a=-\gamma/\alpha\in(-1,0), b=1-1/\alpha\in (0,1/2]$, which is an extension of fractional Brownian motion with Hurst parameter $H\in (1/2, 3/4]$ (corresponding to $\gamma=0$; the process $\xi$ with maximal ranges for the values of the parameters $a$ and $b$ is discussed in \cite{BGT5}). Nevertheless, although the process $\xi$ depends on $\gamma$, its dependence exponent \cite{BGT3} is independent of $\gamma$. The cases $\gamma=\alpha<d$ and $\gamma=d=\alpha$ give a new type of limits (with no counterpart in the homogeneous case), namely, centered, constant (and hence continuous) Gaussian processes  on $(0,\infty)$, discontinuous at $0$.

For $\gamma>d$ the measure $\mu$ is finite, and the results are in sharp contrast to those for $\gamma\leq d$. In this case we give the results for a finite measure $\mu$ in general, and for $d\leq \alpha$ they are akin to the famous limit theorem of Darling and Kac \cite{DK} for the occupation time (without centering), and its generalization to path space by Bingham \cite{Bi}.

All the occupation time  limit theorems are formulated in the context of 
${\cal S}'(\erre^d)$-valued processes, where ${\cal S}'(\erre^d)$ is the usual space of tempered distributions (dual of the space ${\cal S}(\erre^d)$ of smooth rapidly decreasing functions). In some cases the limit process is of the form $\lambda$ multiplied by  a real valued process, but in others the limit is  ``truly'' ${\cal S}'(\erre^d)$-valued.
In all cases where the particle motion is recurrent $(d\leq\alpha)$, the spatial structure of the limit process is $\lambda$, independently of $\gamma$.

The methods of proof for the fluctuation limit theorems are analogous to those developed  in \cite{BGT1, BGT2}, with some new technical complexities because the measure $\mu$ is not invariant for the semigroup of the motion. On the other hand, there is a significant difference in the tightness proofs, as
they require estimates for moments of arbitrary high order (whereas in \cite{BGT1, BGT2} order 2 or 4 was enough). For the results of Darling-Kac type we proceed similarly as \cite{Bi}. However, in our setting the  uniform convergence condition (A) for that kind of results is not satisfied, and some additional work is needed.

Convergence in distribution in the space of continuous functions $C([0,\tau], {\cal S}'(\erre^d))$ for any $\tau>0$ is denoted by $\flecha$. In some cases the interval $[0,\tau]$ is replaced by 
$[\varepsilon,\tau], 0<\varepsilon<\tau$, because the limit process is discontinuous at $t=0$.

The duality between the spaces ${\cal S}'(\erre^d)$ and ${\cal S}(\erre^d)$ is denoted by $\langle\cdot,\cdot\rangle$.

Generic constants are written $C,C_1,C_2,\ldots$, with possible dependencies in parenthesis.

In section 2 we present the results. Section 3 begins with an explanation of the general method used for the proofs of the occupation time fluctuation limits, and then we prove most of the results. Some proofs that are similar to others are omitted, with some  comments.

The branching particle systems in the inhomogeneous case 
produce fewer results, but there are other kinds of difficulties related to extinction of the system. These results will be presented elsewhere.
\vglue1cm
\noindent
{\bf 2 Results}
\setcounter{section}{2}
\setcounter{equation}{0}
\vglue.5cm
Let $N$ and $X_T$ be the processes defined in (\ref{eq:1.1}) and (\ref{eq:1.2}). As stated in the Introduction, for simplicity most of our results are formulated for $\mu$ of the form
\begin{equation}
\label{eq:2.1}
\mu(dx)=\frac{dx}{1+|x|^\gamma},\quad \gamma>0.
\end{equation}
Note that $\mu$ is finite for $d<\gamma$. More general measures $\mu$ will be considered later in this section. In the theorems below, $K$ is a number depending on $\alpha,d,\mu$, which may vary from case to case, and may be computed explicitly in each specific case.

Different arrangements of  $\alpha,\gamma,d$  yield different results, and we order them according to the relationship between $\gamma$ and $d$. We start with $\gamma<d$.

\begin{theorem}
\label{t:2.1}
Let $\gamma<d<\alpha$ (hence $d=1$) and 
\begin{equation}
\label{eq:2.2}
F_T=T^{1-(d+\gamma)/2\alpha}.
\end{equation}
Then $X_T\,\,\flecha\,\, K\lambda\xi$ as $T\to\infty$, where $\xi$ is a real centered Gaussian process with covariance
\begin{equation}
\label{eq:2.3}
E\xi_t\xi_s=\int^{t\wedge s}_0u^{-\gamma/\alpha}\left((t-u)^{1-d/\alpha}+(s-u)^{1-d/\alpha}\right)du.
\end{equation}
\end{theorem}

This theorem is a generalization of Theorem 2.1 in \cite{BGT1}, which corresponds to  $\gamma=0$.
\vglue-.2cm
\begin{remark}
\label{r:2.2}
{\rm The following properties of the process $\xi$ are easy to obtain.

(a) For $\gamma=0, \xi$ is a fractional Brownian motion with Hurst parameter $1-1/2\alpha$ \cite{BGT1}.

(b)$\xi$ is self-similar with index $1-(1+\gamma )/2\alpha$. This is immediate from (2.3), but more generally, from the form of the fluctuation process given by (1.2), it follows that if $F_T$ has the form $T^\kappa f(T)$, where $f$ is a function slowly varying at infinity and $\kappa\geq 0$ (as it is in our cases), then the limit process is self-similar with index $\kappa$ .

(c) $\xi$ is a long range dependence process where
$$E(\xi_{s+T}-\xi_{t+T})(\xi_r-\xi_v)=O(T^{-1/\alpha})\quad{\rm as}\quad T\to\infty,$$
for $0\leq r<v, 0\leq s<t$ .
Note that the dependence exponent \cite{BGT3} $1/\alpha$ is independent of $\gamma$.

(d) $\xi$ is not a Markov process and not a semimartingale. The non-semimartingale property can be proved by Lemma 2.1 in \cite{BGT0}.}
\end{remark}

The next two theorems are generalizations of Theorem 2.1 in \cite{BGT2} (for $\gamma=0$).

\begin{theorem}
\label{t:2.3}
Let $\gamma<d=\alpha\quad (=1$ or $2)$  and
\begin{equation}
\label{eq:2.4}
F_T=(T\log T)^{1/2}T^{-\gamma/2\alpha}.
\end{equation}
Then $X_T\,\,\flecha\,\, K\lambda\beta$ as $T\to\infty$, where $\beta $ is an inhomogeneous real Wiener process with covariance
\begin{equation}
\label{eq:2.5}
E\beta_t\beta_s=\frac{(t\wedge s)^{1-\gamma/\alpha}}{1-\gamma/\alpha}.
\end{equation}
\end{theorem}

\begin{theorem}
\label{t:2.4}
Let $\gamma<\alpha<d$ and
\begin{equation}
\label{eq:2.6}
F_T=T^{(1-\gamma/\alpha)/2}.
\end{equation}
Then $X_T\,\,\flecha\,\, KW$ as $T\to\infty$, where $W$ is an ${\cal S}'(\erre^d)$-valued time inhomogeneous Wiener process with covariance functional
\begin{equation}
\label{eq:2.7}
E\langle W(t),\varphi_1\rangle\langle W(s),\varphi_2\rangle=\frac{(t\wedge s)^{1-\gamma/\alpha}}{1-\gamma/\alpha}\int_{\errita^d}
\varphi_1(x)G\varphi_2(x)dx, \quad \varphi_1,\varphi_2\in{\cal S}(\erre^d),
\end{equation}
where $G$ is the $\alpha$-potential operator, i.e.
\begin{equation}
\label{eq:2.8}
G\varphi(x)=C_{\alpha,d}\int_{\errita^d}\frac{\varphi(y)}{|x-y|^{d-\alpha}}dy,
\end{equation}
with $C_{\alpha,d}=\Gamma(\frac{d-\alpha}{2})(2^\alpha\pi^{\frac{d}{2}}\Gamma(\frac{\alpha}{2}))^{-1}.$
\end{theorem}

The analogy with the case $\gamma=0$ breaks down for ``large'' $\gamma$, i.e., $\gamma\geq \alpha$.

\begin{theorem}
\label{t:2.5}
Let $\gamma=\alpha<d$ and
\begin{equation}
\label{eq:2.9}
F_T=(\log T)^{1/2}.
\end{equation}
Then $X_T\Rightarrow KX$ in $C([\varepsilon,\tau],{\cal S}'(\erre^d))$ as $T\to\infty$ for any $0<\varepsilon<\tau$, where $X$ is an ${\cal S}'(\erre^d)$-valued Gaussian process constant in time on $(0,\infty), X(t)\equiv X(1)$, and $X(1)$ is centered with covariance functional
\begin{equation}
\label{eq:2.10}
E\langle X(1),\varphi_1\rangle\langle X(1),\varphi_2\rangle=\int_{\errita^d}
\varphi_1(x)G\varphi_2(x)dx,\quad \varphi_1,\varphi_2\in{\cal S}(\erre^d),
\end{equation}
where $G$ is given by (\ref{eq:2.8}).
\end{theorem}

Note that the limit process is discontinuous at $t=0$ since $X_T(0)=0$.

To complete the case $\gamma<d$ it remains to consider $\alpha<\gamma<d$. It turns out, however, that if $\alpha<d$ and $\alpha<\gamma$, then the relationship between $\gamma$ and $d$ is irrelevant. In this case the total occupation time
is bounded, so it does not make sense to investigate the fluctuation process. More precisely, we have the following simple proposition.

\begin{proposition}
\label{p:2.6}
Let $\alpha<d$ and $\alpha<\gamma$. Then
$$E\int^\infty_0\langle N_s,\varphi\rangle ds<\infty,\quad\varphi\in{\cal S}(\erre^d).$$
\end{proposition}

We now proceed to the critical case $\gamma=d$.

\begin{theorem}
\label{t:2.7}
Let $1=d=\gamma<\alpha$ and
\begin{equation}
\label{eq:2.11}
F_T=T^{1-d/\alpha}(\log T)^{1/2}.
\end{equation}
Then $X_T\,\,\flecha\,\, K\lambda\xi$ as $T\to\infty$, where $\xi $ is as in Theorem 2.1.
\end{theorem}

The next case is ``doubly critical''.

\begin{theorem}
\label{t:2.8}
Let $\gamma=d=\alpha\quad (=1$ or $2)$ and
\begin{equation}
\label{eq:2.12}
F_T=(\log T)^{3/2}.
\end{equation}
Then $X_T\Rightarrow K\lambda\eta$ in $C([\varepsilon,\tau], {\cal S}'(\erre^d))$ as $T\to\infty$ for any $0<\varepsilon<\tau$, where $\eta$ is a real Gaussian process constant in time on $(0,\infty), \eta_t\equiv\eta_1$, and 
$\eta_1$ is standard normal.
\end{theorem}

Here,  as in Theorem 2.5, the limit process is discontinous at $t=0$.

So far we were assuming that $\mu$ is of the form (\ref{eq:2.1}). It is rather clear that for $\gamma<d$ we can take $\mu(dx)=|x|^{-\gamma} dx$. Moreover, a careful analysis of the proofs shows that in this case  $\mu$ can have a more general form, given in the following proposition. (For the case $\gamma=d$, see the discussion after the proof of Proposition 2.9).

\begin{proposition}
\label{p:2.9}
All the previous results for the case $\gamma<d$ remain true (with possibly different constants $K$) for an intensity measure $\mu$ of the form
\begin{equation}
\label{eq:2.13}
\mu(dx)=\nu(dx)+\frac{h(x)}{1+|x|^\gamma}dx,
\end{equation}
where $\nu$ is a finite measure, and $h$ is a nonnegative bounded function such that there exists a strictly positive limit
\begin{equation}
\label{eq:2.14}
\lim_{R\to\infty}\frac{1}{R^d}\int_{|x|\leq R}h(x)dx.
\end{equation}
\end{proposition}

It seems interesting and perhaps unexpected that it is not sufficient to assume that $\mu(dx)=g(x)dx$ with
\begin{equation}
\label{eq:2.15}
\frac{C_1}{1+|x|^\gamma}\leq g(x)\leq \frac{C_2}{1+|x|^\gamma}.
\end{equation}
We have the following counterexample.

\begin{ex}
\label{e:2.10}{\rm 
Let $\gamma<d<\alpha \,\, (d=1)$ and let $\mu$ be of the form (\ref{eq:2.13}) with $\nu\equiv 0$, and
$$
h(x)=\left\{\begin{array}{lll}
1&{\rm for}&|x|\leq 4,\\
1&{\rm for}&(2k)^{2k}<|x|\leq (2k+1)^{2k+1},\\
2&{\rm for}&(2k+1)^{2k+1}<|x|\leq(2(k+1))^{2(k+1)},
\end{array}\right.
$$%
$k=1,2,\ldots$. The limit (\ref{eq:2.14}) does not exist for this measure, whereas (\ref{eq:2.15}) obviously holds. The only nontrivial normalization (cf. Theorem 2.1) is that given by (\ref{eq:2.2}), but we will explain later that the corresponding $X_T$ does not converge as $T\to\infty$.}
\end{ex}

There remains the case $\gamma>d$. Here the situation changes dramatically and the results are of an entirely different nature. In particular, they do not depend on $\gamma$, but only on the fact that the measure $\mu$ is finite. Therefore we will formulate our results for a general finite measure $\mu$. It turns out that the appropriate normalization is
\begin{equation}
\label{eq:2.16}
F_T=T^{1-1/\alpha}
\end{equation}
if $1=d<\alpha$, and
\begin{equation}
\label{eq:2.17}
F_T=\log T
\end{equation}
if $d=\alpha$.

It is easy to see that in both cases $\frac{1}{F_T}\int^T_0\langle N_s,\varphi\rangle ds$ converges to a finite limit as $T\to\infty,\varphi\in{\cal S}(\erre^d)$, hence there is no reason to consider fluctuation processes and it suffices to investigate the occupation process
\begin{equation}
\label{eq:2.18}
Y_T(t)=\frac{1}{F_T}\int^{Tt}_0N_sds.
\end{equation}

For a given $\alpha>1$, let $L$ denote the local time process (at $0$) of a standard real $\alpha$-stable process. See, e.g., \cite{Be} for properties of $L$. In particular, $L$ is a continuous increasing process, $L(0)=0$.

The relation between the processes $Y_T$ and $L$ is given in the following theorem.

\begin{theorem}
\label{t:2.11}
Let $1=d<\alpha$ and $\mu$ be a finite measure on $\erre$. Let $L_1,L_2,\ldots$ be independent copies of $L$ and let $\nu$ be a Poisson random variable with parameter $\mu(\erre)$ independent of $L_1,L_2,\ldots$. Then for $F_T$ defined by 
(\ref{eq:2.16}),
$$Y_T\,\,\flecha\,\, K\sum_{j\leq \nu}L_j\lambda$$
as $T\to\infty$.
\end{theorem}

This theorem is based on the following lemma, which is of interest by itself.

\begin{lemma}
\label{l:2.12}
Let $\alpha, d$ and $F_T$ be as before, let $\zeta$ be a real standard $\alpha$-stable process, and denote by $Z_T$ its normalized occupation process, i.e.,
\begin{equation}
\label{eq:2.19}
\langle Z_T(t),\varphi\rangle=\frac{1}{F_T}\int^{tT}_0\varphi(\zeta_s)ds,\quad\varphi\in{\cal S}(\erre^d), t\geq 0.
\end{equation}
Then
\begin{equation}
\label{eq:2.20}
Z_T\,\,\flecha\,\, KL\lambda
\end{equation}
as $T\to\infty$.
\end{lemma}

This lemma is closely related to the famous Darling-Kac result \cite{DK}. Their theorem was generalized  by Bingham \cite{Bi}, who obtained the limit in path space for fixed positive $\varphi$ with compact support. Fitzsimmons and Getoor 
\cite{FG1} mention this limit for fixed general $\varphi$. We will present an outline of a proof of the lemma in the next section.

It remains to consider the case $d=\alpha$.

\begin{theorem}
\label{t:2.13}
Let $d=\alpha \quad (=1$ or $2)$ and $\mu$ be a finite measure. Let $\rho_1,\rho_2,\ldots$ be i.i.d. standard exponential random variables and $\nu$  a Poisson random variable with parameter $\mu(\erre)$, independent of $\rho_1,\rho_2,\ldots$. Then for $F_T$ defined by (\ref{eq:2.17}),
$$Y_T\Rightarrow K\sum_{j\leq\nu}\rho_j\lambda$$
in $C([\varepsilon,\tau],{\cal S}'(\erre^d))$ as $T\to\infty$ for any $0<\varepsilon<\tau$.
\end{theorem}

So the limit process is constant in time on $(0,\infty)$.
\vglue1cm
\noindent
{\bf 3 Proofs}
\setcounter{section}{3}
\setcounter{equation}{0}

\subsection{General scheme}

We describe a general method used in the proofs of Theorems 2.1-2.9. 
Given $0<\alpha\leq 2$,
let ${\cal T}_t$ denote the transition semigroup of the standard $\alpha$-stable process $\zeta$ in $\erre^d$, i.e. ${\cal T}_t\varphi=p_t *\varphi$, where $p_t$ is the transition density of $\zeta$ .

It is well known that, by the Poisson property, 
\begin{equation}
\label{eq:3.0}
E\langle N_t,\varphi\rangle=\int_{\errita^d}{\cal T}_t\varphi(x)\mu(dx),\quad \varphi\in{\cal S}(\erre^d).
\end{equation}

For a continuous ${\cal S}'(\erre^d)$-valued process $X$ we define an ${\cal S}'(\erre^{d+1})$ random variable $\widetilde{X}$ by
\begin{equation}
\label{eq:3.1}
\langle\widetilde{X},\Phi\rangle=
\int^\tau_0\langle X(t),\Phi(\cdot, t)\rangle dt,\quad\Phi\in{\cal S}(\erre^{d+1}).
\end{equation}
As explained in \cite{BGR}, in order to prove $X_T\,\,\flecha\,\, X$, where $X$ is the limit process occuring in the specific theorem, it suffices to show that 
\begin{equation}
\label{eq:3.2}
\langle\widetilde{X}_T,\Phi\rangle\Rightarrow \langle\widetilde{X},\Phi\rangle,\quad\Phi\in{\cal S}(\erre^{d+1}),
\end{equation}
and that the family $\{\langle X_T,\varphi\rangle\}_{T\geq 2}$ is tight in 
$C([0,\tau],\erre)$ for any $\tau>0$, for each $\varphi\in{\cal S}(\erre^d)$.

This scheme should be modified in an obvious way if we consider $C([\varepsilon,\tau],{\cal S}'(\erre^d))$. Without loss of generality we will always assume $\tau=1$.

Since the limits are Gaussian, in order to obtain 
(\ref{eq:3.2}) it suffices to show that
\begin{equation}
\label{eq:3.3}
\lim_{T\to\infty}Ee^{-\langle\widetilde{X}_T,\Phi\rangle}=Ee^{-\langle\widetilde{X},\Phi\rangle}
\end{equation}
for any nonnegative $\Phi\in{\cal S}(\erre^{d+1})$ (see, e.g., \cite{BGT1}).

Given such $\Phi$ we denote
\begin{equation}
\label{eq:3.4}
\Psi(x,t)=\int^1_t\Phi(x,s)ds,\quad\Psi_T(x,t)=\frac{1}{F_T}\Psi\biggl(x,\frac{t}{T}\biggr).
\end{equation}

By (\ref{eq:1.2}), (\ref{eq:3.0}), (\ref{eq:3.1}),  we have
\begin{equation}
\label{eq:3.5}
\langle\widetilde{X}_T,\Phi\rangle=\int^T_0\langle N_u,\Psi_T(\cdot,u)\rangle du-\int^T_0\int_{\errita^d}{\cal T}_u\Psi_T(\cdot,u)\mu(dx)du.
\end{equation}
Hence,  by the Poisson property,
\begin{equation}
\label{eq:3.6}
Ee^{-\langle\widetilde{X}_T,\Phi\rangle}={\rm exp}\bigg\{\int^T_0\int_{\errita^d}{\cal T}_u\Psi_T(\cdot,u)(x)\mu(dx)du\biggr\}
{\rm exp}\biggl\{-\int_{\errita^d}v_T(x,T)\mu(dx)\biggr\},
\end{equation}
where
\begin{equation}
\label{eq:3.7}
v_T(x,t)=1-E{\rm exp}\left\{
-\int^t_0\Psi_T(x+\zeta_u,T-t+u)du\right\},\quad 0\leq t\leq T.
\end{equation}
(Recall that $\zeta$ is the standard $\alpha$-stable process).

We know that (repeating the argument in \cite{BGT1} for $V=0$), by the 
Feynman-Kac formula, $v_T$ satisfies
\begin{equation}
\label{eq:3.8}
v_T(x,t)=\int^t_0{\cal T}_{t-s}(\Psi_T(\cdot,T-s)(1-v_T(\cdot,s)))(x)ds.
\end{equation}

We will often use an immediate consequence of (\ref{eq:3.7}) and 
(\ref{eq:3.8}):
\begin{equation}
\label{eq:3.9}
v_T(x,t)\leq\int^t_0{\cal T}_{t-s}\Psi_T(\cdot,T-s)(x)ds.
\end{equation}

Putting (\ref{eq:3.8}) into (\ref{eq:3.6}) and then using (\ref{eq:3.8}) once more we obtain
\begin{equation}
\label{eq:3.10}
Ee^{-\langle\widetilde{X}_T,\Phi\rangle}=e^{I(T)-I\!I(T)},
\end{equation}
where
\begin{equation}
\label{eq:3.11}
I(T)=\int_{\errita^d}\int^T_0{\cal T}_{T-s}\left(\Psi_T(\cdot,T-s)\int^s_0{\cal T}_{s-u}\Psi_T(\cdot,T-u)du\right)(x)ds\mu(dx)
\end{equation}
and
\begin{equation}
\label{eq:3.12}
I\!I(T)=\int_{\errita^d}\int^T_0{\cal T}_{T-s}\left(\Psi_T(\cdot,T-s)\int^s_0{\cal T}_{s-u}(\Psi_T(\cdot,T-u)v_T(\cdot,u))du\right)(x)ds\mu(dx).
\end{equation}

To prove (\ref{eq:3.3}) we will show that
\begin{equation}
\label{eq:3.13}
\lim_{T\to\infty}e^{I(T)}=Ee^{-\langle\widetilde{X},\Phi\rangle},
\end{equation}
and
\begin{equation}
\label{eq:3.14}
\lim_{T\to\infty}I\!I(T)=0.
\end{equation}

For simplicity we will prove (\ref{eq:3.13}) and (\ref{eq:3.14}) for $\Phi$ of the form

\begin{equation}
\label{eq:3.15}
\Phi(x,t)=\varphi(x)\psi(t),\quad \varphi\in {\cal S}(\erre^d),\psi\in{\cal S}(\erre),\quad\varphi,\psi\geq 0.
\end{equation}
For such $\Phi$ it will be convenient to denote
\begin{equation}
\label{eq:3.16}
\chi(t)=\int^1_t\psi(s)ds,\quad\chi_T(t)=\chi\biggl(\frac{t}{T}\biggr),
\end{equation}
then
\begin{equation}
\label{eq:3.17}
\Psi_T(x,t)=\frac{1}{F_T}\varphi(x)\chi_T(t).
\end{equation}

Note that expressions $I(T)$ and $I\!I(T)$ have more complicated forms than those corresponding to $\gamma=0$ \cite{BGT1, BGT2}, since the measure $\mu$ is not invariant under ${\cal T}_t$ and it is infinite, so in particular the Fourier transform technique we have used before is not applicable.

In order to prove tightness of $\{\langle X_T,\varphi\rangle\}_{T\geq 2}$ for a given $\varphi\in{\cal S}(\erre^d), \varphi\geq 0$ (it suffices to take $\varphi$ nonnegative), we need a formula for the Laplace transform of $\langle X_T(t_2)-X_T(t_1),\varphi\rangle$ for $0\leq t_1<t_2\leq 1$. We take $\Psi_{T,n}$ of the form  (\ref{eq:3.17}) with $\theta\varphi$ instead of $\varphi\quad (\theta>0)$, and with smooth $\chi_n$ approximating $\chi=\UNO_{[t_1,t_2]}$. Using 
 (\ref{eq:3.8}),  (\ref{eq:3.6}) and  (\ref{eq:3.5}) and letting $n\to\infty$ we obtain
\begin{equation}
\label{eq:3.18}
Ee^{-\theta\langle X_T(t_2)-X_T(t_1),\varphi\rangle}=e^{H_T(\theta)},
\end{equation}
where
\begin{equation}
\label{eq:3.19}
H_T(\theta)=\frac{\theta}{F_T}\int_{\errita^d}\int^T_0{\cal T}_{T-s}(\varphi\chi_T(t-s)v_{\theta,T}(\cdot,s))(x)ds\mu(dx),
\end{equation}
and $v_{\theta,T}$ is defined by  (\ref{eq:3.7}) for $\Psi_T(x,t)=\theta\varphi(x)\chi_T(t)$.
This $v_{\theta,T}$ also satisfies  (\ref{eq:3.8}).

Unlike \cite{BGT2}, where fourth moments were employed, we need moments of $\langle X_T(t_2)-X_T(t_1),\varphi\rangle$ of arbitrary high order. By  (\ref{eq:3.18}) we have
\begin{equation}
\label{eq:3.20}
E\langle X_T(t_2)-X_T(t_1),\varphi\rangle^k=(-1)^k\frac{d^k}{d\theta^k}e^{H_T(\theta)}\biggl|_{\theta=0},\quad k=1,2,\ldots.
\end{equation}
Using  (\ref{eq:3.8}) and  (\ref{eq:3.19}) we have
\begin{eqnarray}
H_T(0)&=&0,\nonumber\\
H'_T(0)&=&0,\nonumber\\
\label{eq:3.21}
H^{(k)}_T(0)&=&(-1)^k\frac{k!}{F^k_T}\int_{\errita^d}\int^T_0\int^{s_k}_0\ldots\int^{s_2}_0{\cal T}_{T-s_k}(\varphi{\cal T}_{s_k-s_{k-1}}(\varphi\ldots{\cal T}_{s_2-s_1})\ldots)(x)\nonumber\\
&&\qquad \times \chi_T(T-s_k)\ldots\chi_T(T-s_1)ds_1\ldots ds_k\mu(dx), \quad k\geq 2.
\end{eqnarray}
By  (\ref{eq:3.20}) and  (\ref{eq:3.21}), tightness will be proved if we show that there exists $\delta>0$ such that
\begin{equation}
\label{eq:3.22}
|H^{(k)}_T(0)|\leq C(k,\varphi)(t_2-t_1)^{k\delta}\quad{\rm for}\quad k=2,3,\ldots.
\end{equation}

The scheme described above is employed in the proofs of all results for $\gamma\leq d$. The proof of each specific case, however, requires slightly different and non-trivial calculations; nevertheless, for brevity we will omit some proofs, concentrating on arguments which are either the most typical or the most involved.

\subsection{Proof of Theorem 2.1}

We will prove the theorem for $\mu(dx)=|x|^{-\gamma}dx$ since in this case the formulas are slightly simpler. It will be obvious that the same type of argument applies for $\mu$ of the form (2.1). It is easy to see that in this case the right hand side of  (\ref{eq:3.13}) with $\Phi$ given by  (\ref{eq:3.15}) is of the form 
\begin{equation}
\label{eq:3.23}
{\rm exp}\biggl\{K_1\biggl(\int\varphi(x)dx\biggr)^2\int^1_0\int^u_0(u-s)^{-d/\alpha}s^{-\gamma/\alpha}\chi(s)\chi(u)dsdu\biggr\}.
\end{equation}

Using  (\ref{eq:3.11}),  (\ref{eq:3.17}) and substituting $u'=1-u/T,s'=1-s/T$ we obtain
$$I(T)=\frac{T^2}{F^2_T}\int_{\errita^d}\int^1_0\int^u_0\int_{\errita^{2d}}p_{Ts}(x-y)\varphi(y)p_{T(u-s)}(y-z)\varphi(z)\chi(s)\chi(u)|x|^{-\gamma}dydzdsdudx.$$
We apply the self similarity of the $\alpha$-stable density, i.e.
\begin{equation}
\label{eq:3.24}
p_{at}(x)=a^{-d/\alpha}p_t(xa^{-1/\alpha}),
\end{equation}
substitute $x'=xT^{-1/\alpha}$ and use (2.2), then
\begin{equation}
\label{eq:3.25}
I(T)=\int_{\errita^d}g_T(x)|x|^{-\gamma}dx,
\end{equation}
where
\begin{equation}
\label{eq:3.26}
g_T(x)=\int^1_0\int^u_0\int_{\errita^{2d}}\chi(s)\chi(u)p_s(x-yT^{-1/\alpha})p_{u-s}((y-z)T^{-1/\alpha})\varphi(y)\varphi(z)dydzdsdu.
\end{equation}
By self similarity again the integrand in (\ref{eq:3.26}) is bounded by $Cs^{-d/\alpha}(u-s)^{-d/\alpha}\varphi(y)\varphi(z)$, which is integrable since $d<\alpha$. Hence, by the dominated convergence theorem we obtain
\begin{eqnarray}
\lefteqn{\lim_{T\to\infty}g_T(x)=g_\infty(x)=\int^1_0\int^u_0\int_{\errita^{2d}}\chi(s)\chi(u)p_s(x)p_{u-s}(0)\varphi(y)\varphi(z)dydzdsdu}\nonumber\\
\label{eq:3.27}
&&=p_1(0)\int^1_0\int^u_0\chi(s)\chi(u)s^{-d/\alpha}(u-s)^{-d/\alpha}p_1(xs^{-1/\alpha})dsdu\left(\int_{\errita^d}\varphi(y)dy\right)^2.
\end{eqnarray}
It is well known that
\begin{equation}
\label{eq:3.28}
p_1(x)\leq \frac{C}{1+|x|^{d+\alpha}},
\end{equation}
hence, by (\ref{eq:3.24}) we easily deduce that for $d<\alpha$,
\begin{equation}
\label{eq:3.29}
\int^1_0p_s(x)ds\leq \frac{C_1}{1+|x|^{d+\alpha}}.
\end{equation}
This and an obvious estimate
\begin{equation}
\label{eq:3.30}
\frac{1}{1+|x-w|^{d+\alpha}}\leq C_2\frac{1+|w|^{d+\alpha}}{1+|x|^{d+\alpha}}
\end{equation}
(for $w=yT^{-1/\alpha}$) imply that
\begin{equation}
\label{eq:3.31}
g_T(x)\leq C_3\frac{1}{1+|x|^{d+\alpha}}.
\end{equation}
By (\ref{eq:3.25}), (\ref{eq:3.27}), (\ref{eq:3.31}) and taking into account (\ref{eq:3.23}), we obtain (\ref{eq:3.13}).

This completes the proof of (\ref{eq:3.2}).

We proceed to the proof of (\ref{eq:3.14}). We use (\ref{eq:3.12}), 
(\ref{eq:3.9}), (\ref{eq:3.17}) and boundedness of $\chi$ to get
\begin{equation}
\label{eq:3.32}
I\!I(T)\leq\frac{C}{F^3_T}\int_{\errita^d}\int^T_0\int^s_0\int^u_0{\cal T}_{T-s}(\varphi{\cal T}_{s-u}(\varphi{\cal T}_{u-v}\varphi))(x)dvduds|x|^{-\gamma}dx.
\end{equation}
We substitute $v'=\frac{u-v}{T}$, then $u'=\frac{s-u}{T}$, then 
$s'=\frac{1-s}{T}$, and we increase the time intervals to $[0,1]$, obtaining
\begin{eqnarray}
I\!I(T)&\leq& C\frac{T^3}{F^3_T}
\int_{\errita^d}\int_{\errita^{3d}}\int^1_0p_{Ts}(x-y)ds
\varphi(y)\int^1_0p_{Tu}(y-z)du\varphi(z)\int^1_0p_{Tv}(z-w)dv\nonumber\\
\label{eq:3.33}
&&\times \varphi(w)dydzdw|x|^{-\gamma}dx.
\end{eqnarray}
Denote

\begin{equation}
\label{eq:3.34}
f(x)=\int^1_0p_s(x)ds
\end{equation}
and
\begin{equation}
\label{eq:3.35}
\widetilde{\varphi}_T(x)=T^{d/\alpha}\varphi(T^{1/\alpha}x).
\end{equation}
Note that $f$ is integrable, and by (\ref{eq:3.29}) it is bounded, and
\begin{equation}
\label{eq:3.36}
\int_{\errita^d}\widetilde{\varphi}_T(x)dx=\int_{\errita^d}\varphi(x)dx.
\end{equation}

Using (\ref{eq:3.24}), substituting 
$x'=xT^{-1/\alpha},\,y'=yT^{-1/\alpha},
z'=zT^{-1/\alpha},w'=wT^{-1/\alpha}$,
we write (\ref{eq:3.33}) as
\begin{equation}
\label{eq:3.37}
I\!I(T)\leq C
\frac{T^{3-(2d+\gamma)/\alpha}}{F^3_T}
\int_{\errita} a_T(x)|x|^{-\gamma}dx,
\end{equation}
where
$$a_T(x)=f*(\widetilde{\varphi}_T*(f*(\widetilde{\varphi}_T(f*\widetilde{\varphi}_T))))(x).$$

The properties of $f$ and $\widetilde{\varphi}_T$ easily imply that
$$\sup_T\int_{\errita^d}a_T(x)dx<\infty\quad{\rm and}\quad 
\sup_T\sup_{x\in\errita^d}a_T(x)<\infty.$$
Hence (\ref{eq:3.14}) follows from (\ref{eq:3.37}) since 
$\gamma<d$ and $T^{3-(2d+\gamma)/\alpha}/
F^3_T\to 0$ (see (\ref{eq:2.2})).

According to the general scheme, in order to prove tightness we show (\ref{eq:3.22}). We substitute $s'_j=1-\frac{s_j}{T}$ in (\ref{eq:3.21}) and we obtain
\begin{eqnarray}
\label{eq:3.38}
\lefteqn{|H^{(k)}(0)|}\\
&=&k!\frac{T^k}{F^k_T}\int_{\errita^d}\int^1_0\int^1_{s_k}
\ldots\int^1_{s_2}{\cal T}_{Ts_k}(\varphi{\cal T}_{T(s_{k-1}-s_k)}(\varphi\ldots)\ldots)(x)\chi(s_k)\ldots\chi g(s_1)
ds_1\ldots ds_k\frac{1}{|x|^\gamma}dx.\nonumber
\end{eqnarray}
We need the following estimate:
\begin{equation}
\label{eq:3.39}
\frac{T}{F_T}\int^1_s{\cal T}_{T(u-s)}\varphi(y)\chi(u)du\leq 
CT^{(\gamma -d)/2\alpha}(t_2-t_1)^{1-d/\alpha}
\end{equation}
(recall that $\chi=\UNO_{[t_1,t_2]}$). By 
(\ref{eq:2.2}), (\ref{eq:3.24}) and boundedness of $p_1$ we have
\begin{equation}
\label{eq:3.40}
\frac{T}{F_T}\int^1_s{\cal T}_{T(u-s)}\varphi(y)\chi(u)du\leq 
CT^{(\gamma-d)/2\alpha}\int_{\errita^d}\varphi(z)dz\int^1_s(u-s)^{-d/\alpha}\chi(u)du.
\end{equation}
Hence (\ref{eq:3.39}) follows. We  iterate (\ref{eq:3.39}) $k-1$ times in 
(\ref{eq:3.38}), estimate $(T^{(\gamma -d)/2\alpha})^{k-2}$ by $1$, arriving at
\begin{eqnarray*}
|H^{(k)}(0)|&\leq& C(t_2-t_1)^{(k-1)(1-d/\alpha)}T^{(\gamma -d)/2\alpha}\frac{T}{F_T}\int_{\errita^d}\int^1_0\int_{\errita^d}p_{Ts_k}(x-y)\varphi(y)dyds_k|x|^{-\gamma}dx\\
&=&C(t_2-t_1)^{(k-1)(1-d/\alpha)}\int_{\errita^d}f*\widetilde{\varphi}_T(x)|x|^{-\gamma}dx,
\end{eqnarray*}
where we have used (\ref{eq:2.2}), self-similarity and the usual substitutions $x'=xT^{-1/\alpha}, y'=yT^{-1/\alpha}$, where $f$ and $\widetilde{\varphi}_T$ are defined by (\ref{eq:3.34}) and (\ref{eq:3.35}). Hence we obtain (\ref{eq:3.22}) by the properties of $f$ and $\widetilde{\varphi}_T$. 

The proof of the theorem is complete. $\hfill\Box$

\subsection{Some properties of the $\alpha$-stable semigroup in the critical case $d=\alpha$.}

We will need the following facts, valid for $d=\alpha,\varphi\in{\cal S}(\erre^d),\varphi\geq 0$:
\begin{eqnarray}
\label{eq:3.41}
&&\sup_{T>2}\sup_{x\in\errita^d}\frac{1}{\log T}\int^T_0{\cal T}_u\varphi(x)du<\infty,\\
\label{eq:3.42}
&&\lim_{T\to\infty}\frac{1}{\log T}\int^T_0{\cal T}_u\varphi(x)du=p_1(0)\int_{\errita^d}\varphi(y)dy,\\
\label{eq:3.43}
&&\lim_{T\to\infty}\frac{1}{\log T}\int_{|x|^d>T}\int^T_0{\cal T}_u\varphi(x)du|x|^{-d}dx=0.
\end{eqnarray}
These properties are perhaps known but we have not been able to find  references for them, so we show briefly how to derive them.

To prove (\ref{eq:3.41}) and (\ref{eq:3.42}), it is clear that it suffices to consider $\int^T_1$. We use self-similarity, then make a substitution which turns out to be particularly useful in the critical cases and will be applied several times. Namely, we put
\begin{equation}
\label{eq:3.44}
u'=\frac{\log u}{\log T},
\end{equation}
obtaining
$$\frac{1}{\log T}\int^T_1{\cal T}_u\varphi(x)du=\int^1_0\int_{\errita^d}p_1\left((x-y)T^{-u/d}\right)\varphi(y)dydu.$$
Hence (\ref{eq:3.41}) and (\ref{eq:3.42}) follow immediately.

To prove (\ref{eq:3.43}) we again replace $\int^T_0$ by $\int^T_1$ and make the substitution (\ref{eq:3.44}). We then have
\begin{eqnarray*}
\lefteqn{\frac{1}{\log T}\int_{|x|^d>T}\int^T_1 {\cal T}_u\varphi(x)du|x|^{-d}
dx}\\
&=&\int_{|x|^d>T}\int^1_0\int_{\errita^d}p_1((x-y)T^{-u/d})\varphi(y)|x|^{-d}dydudx,\\
&\leq& C\int_{|x|^d>T}\int^1_0\int_{\errita^d}\frac{1}{1+|x|^{2d}T^{-2u}}(1+|y|^{2d}T^{-2u})\varphi(y)|x|^{-d}dydudx,
\end{eqnarray*}
where the last estimate follows from (\ref{eq:3.28}) and (\ref{eq:3.30}) (recall that $d=\alpha$). As $\varphi\in{\cal S}(\erre^d)$, this expression is estimated by
$$C_1\int_{|x|^d>T}|x|^{-3d}dx\int^1_0T^{2u}du\leq \frac{C_2}{\log T}$$
by calculus. This proves (\ref{eq:3.43}).

\subsection{Proof of Theorem 2.3}

We will present only an outline of the proof.

Following the general scheme, and again taking for simplicity $\mu(dx)=|x|^{-\gamma}dx$, we prove that (see (\ref{eq:3.11}), (\ref{eq:3.16}),(\ref{eq:3.17}))
\begin{equation}
\label{eq:3.45}
\lim_{T\to\infty}I(T)=K_1\int^1_0s^{-\gamma /\alpha}
\chi^2(s)ds\biggl(\int_{\errita^d}\varphi(x)dx\biggr)^2.
\end{equation}
In (\ref{eq:3.11}) we substitute $u'=T-u, s'=1-\frac{s}{T}$, use self-similarity and put $x'=xT^{-1/\alpha}s^{-1/\alpha}$. By (\ref{eq:2.4}) we obtain
\begin{eqnarray}
\label{eq:3.46}
I(T)&=&\frac{1}{\log T}\int_{\errita^d}\int^1_0\int_{\errita^d}\int^{T(1-s)}_0\chi(s)\chi\biggl(s+\frac{u}{T}\biggr)s^{-\gamma /\alpha }p_1(x-y
s^{-1/\alpha}T^{-1/\alpha})\\
&& \times \varphi(y){\cal T}_u\varphi(y)|x|^{-\gamma}dudydsdx.\nonumber
\end{eqnarray}
Using
\begin{equation}
\label{eq:3.47}
\sup_{z\in\errita^d}\int_{\errita^d}p_1(x+z)|x|^{-\gamma}dx<\infty,
\end{equation}
it is easy to see that
$$\lim_{T\to\infty}I(T)=\lim_{T\to\infty}I'(T),$$
where
\begin{eqnarray*}
I'(T)&=&\frac{1}{\log T}
\int_{\errita^d}\int^{1-\frac{1}{T}}_0
\int_{\errita^d}\int^{T(1-s)}_1\chi(s)\chi\biggl(s+\frac{u}{T}\biggr)
s^{-\gamma/\alpha}p_1(x-ys^{-1/\alpha}
T^{-1/\alpha})\\
&&\times\varphi(y){\cal T}_u\varphi(y)|x|^{-\gamma}dudydsdx.
\end{eqnarray*}

We use self-similarity again and make substitution (\ref{eq:3.44}). Then
\begin{eqnarray*}
I'(T)&=&\int_{\errita^d}\int^{1-\frac{1}{T}}_0\int_{\errita^d}
\int^{\frac{\log T(1-s)}{\log T}}_0\chi(s)\chi(s+T^{u-1})s^{-\alpha/
\gamma}p_1(x-y
s^{-1/\alpha}T^{-1/\alpha})\\
&&\times\varphi(y)\int_{\errita^d}p_1((y-z)T^{-u/\alpha})\varphi(z)dzdudyds|x|^{-\gamma}dx.
\end{eqnarray*}

Now it is clear that the limit of $I'(T)$ should have the form (\ref{eq:3.45}). We omit details.

We proceed to (\ref{eq:3.14}). We use (\ref{eq:3.32}), substitute $v'=u-v$, then $u'=s-u$, then $s'=\frac{1-s}{T}$, and increase the time intervals appropriately, obtaining
\begin{equation}
\label{eq:3.48}
I\!I(T)\leq C\frac{T}{F^3_T}\int_{\errita^d}\int^1_0{\cal T}_{Ts}\left(\varphi\int^T_0{\cal T}_u\biggl(\varphi\int^T_0{\cal T}_v\varphi dv\biggr)du\right)(x)ds|x|^{-\gamma}dx.
\end{equation}
By (\ref{eq:3.41}) applied twice we have
$$I\!I(T)\leq C_1\frac{T(\log T)^2}{F^3_T}\int_{\errita^d}\int^1_0{\cal T}_{Ts}\varphi(x)ds|x|^{-\gamma}dx.$$
Hence (\ref{eq:3.14}) follows by self-similarity, substitution $x'=x(Ts)^{-1/\alpha}$, (\ref{eq:3.47}) and (\ref{eq:2.4}).

Tightness is proved similarly as in Theorem 2.1. Here are the main steps. Instead of (\ref{eq:3.37}) we show that
$$\frac{T}{F_T}\int^1_s{\cal T}_{T(u-s)}\varphi(y)\chi(u)du\leq C(t_2-t_1)^{\frac{1}{4}(1-\gamma/\alpha)},$$
we iterate this estimate $k-2$ times in (\ref{eq:3.40}), and we obtain
$$|H^{(k)}(0)|\leq C_1(t_2-t_1)^{(1-\gamma/\alpha)(k-2)/4}H''(0).$$

Using 
(\ref{eq:3.41}) it is not difficult to prove that
$$H''(0)\leq C_2(t_2^{1-\gamma/\alpha}-t_1^{1-\gamma/\alpha})\leq C_2(t_2-t_1)^{1-\gamma/\alpha}.$$
Hence (\ref{eq:3.22}) follows.  $\hfill\Box$

\subsection{Proof of Theorem 2.4}

Again we give only a sketch of the proof. We make the same substitutions as in the beginning of the previous proof, and by 
(\ref{eq:2.6}) we obtain
$$
I(T)=\int_{\errita^d}\int^1_0\int_{\errita^d}\int^{T(1-s)}_0\chi(s)\chi\biggl(s+\frac{u}{T}\biggr)s^{-\gamma/\alpha}p_1(x-ys^{-1/\alpha}T^{-1/\alpha})\varphi(y){\cal T}_u\varphi(y)|x|^{-\gamma}dudydsdx
$$
(cf. (\ref{eq:3.46})). Hence it is not hard to see that
$$\lim_{T\to\infty}I(T)=\int_{\errita^d}p_1(x)|x|^{-\gamma}dx\int^1_0
s^{-\gamma/\alpha}\chi^2(s)ds\int_{\errita^d}\varphi(y)G\varphi(y)dy,$$
since $\int^\infty_0{\cal T}_u\varphi du=G \varphi$. This implies 
(\ref{eq:3.13}). Next, by (\ref{eq:3.48}) (which is always valid),
$$I\!I(T)\leq C\frac{T}{F^3_T}\int_{\errita^d}\int^1_0{\cal T}_{Ts}(\varphi G(\varphi G\varphi))(x)ds|x|^{-\gamma}dx.$$
Hence (\ref{eq:3.14}) follows by the usual argument since $G\varphi$ is bounded.

Tightness can be proved in the same manner, even easier, as in Theorem 2.3.
$\hfill\Box$

\subsection{Comments on the proofs of Theorems 2.5-2.7}

The proof of Theorem 2.5, though by no means straightforward, is slightly simpler than the proof for the doubly critical case (Theorem 2.8), which will be given in detail. Therefore we omit it. Proposition 2.6 is obtained immediately from 
(\ref{eq:3.0}) and the following estimate, valid for $\varphi\geq 0$:
$$E\int^\infty_0\langle N_s,\varphi\rangle ds=\int_{\errita^d}G\varphi(x)\frac{1}{1+|x|^\gamma}dx\leq C\int_{\errita^d}\frac{1}{1+|x|^{d-\alpha}}\frac{1}{1+|x|^\gamma}<\infty,$$
where we have used the fact that $\sup_x(1+|x|^{d-\alpha})|G\varphi(x)|<\infty$ \cite{I} (Lemma 5.3).

The proof of Theorem 2.7 is similar (but not identical) to the argument carried out for Theorem 2.1. We omit it for brevity.

\subsection{Proof of Theorem 2.8}

We apply the general scheme. By (\ref{eq:3.11}), (\ref{eq:3.15})-(\ref{eq:3.17}) and the substitutions $u'=s-u,$ then $s'=T-s$, we have
\begin{equation}
\label{eq:3.49}
I(T)=I_1(T)+I_2(T)+I_3(T)+I_4(T),
\end{equation}
where
\begin{eqnarray}\label{eq:3.50}
I_1(T)&=&\frac{1}{F^2_T}\int_{1\leq|x|^d\leq T}\int^{T-1}_1\int^{T-s}_1{\cal T}_s(\varphi {\cal T}_u\varphi)(x)\chi\biggl(\frac{s}{T}\biggr)\chi\biggl(\frac{s}{T}+\frac{u}{T}\biggr)\frac{1}{1+|x|^d}dudsdx,\\
\label{eq:3.51}
I_2(T)&=&\frac{1}{F^2_T}\int_{1\leq |x|^d\leq T}\left(\int^T_0\int^{T-s}_0-\int^{T-1}_1\int^{T-s}_1\right)\ldots,\\
\label{eq:3.52}
I_3(T)&=&\frac{1}{F^2_T}\int_{|x|^d>T}\int^T_0\int^{T-s}_0\ldots,\\
\label{eq:3.53}
I_4(T)&=&\frac{1}{F^2_T}\int_{|x|<1}\int^T_0\int^{T-s}_0\ldots,
\end{eqnarray}
where $\ldots$ denotes the same integrand as in $I_1(T)$.

We will show that
\begin{equation}
\label{eq:3.54}
\lim_{T\to\infty}I_1(T)=K_1\chi^2(0)\left(\int_{\errita^d}\varphi(x)dx\right)^2,
\end{equation}
and the remaining integrals converge to $0$.

By (\ref{eq:2.12})
$$I_2(T)\leq\frac{1}{(\log T)^3}\int_{|x|^d\leq T}\left(\int^T_0\int^1_0\ldots duds+\int^1_0\int^T_0\ldots duds\right)dx.$$
Using (\ref{eq:3.41}) we obtain
\begin{equation}
\label{eq:3.55}
I_2(T)\leq\frac{C_1}{(\log T)^2}\int_{|x|^d\leq T}\frac{1}{1+|x|^d}dx\leq \frac{C_2}{\log T}\rightarrow 0 \quad {\rm as}\quad T\rightarrow \infty.
\end{equation}

The fact that
\begin{equation}
\label{eq:3.56}
\lim_{T\to\infty}I_3(T)=0
\end{equation}
follows immediately from (\ref{eq:3.41}) and (\ref{eq:3.43}), and 
\begin{equation}
\label{eq:3.57}
\lim_{T\to\infty}I_4(T)=0
\end{equation}
is also a consequence of (\ref{eq:3.41}).

By (\ref{eq:3.50}), (\ref{eq:2.12}) and (\ref{eq:3.24}) we have
\begin{eqnarray*}
I_1(T)&=&\frac{1}{(\log T)^3}\int_{1\leq|x|^d\leq T}\int_{\errita^{2d}}\int^{T-1}_1\int^{T-s}_1s^{-1}p_1((x-y)s^{-1/d})\varphi(y)u^{-1}p_1((y-z)u^{-1/d})\\
&&\qquad\qquad\qquad \times\varphi(z)\chi\biggl(\frac{s}{T}\biggr)\chi\biggl(\frac{s}{T}+\frac{u}{T}\biggr)\frac{1}{1+|x|^d}dudsdydzdx.
\end{eqnarray*}
We make the substitution (\ref{eq:3.44}) for both $u$ and $s$, obtaining
\begin{eqnarray*}
I_1(T)&=&\frac{1}{\log T}\int_{1\leq|x|^d\leq T}
\int_{\errita^{2d}}\int_0^{\frac{\log(T-1)}{\log T}}
\int_0^{\frac{\log(T-T^s)}{\log T}}
p_1((x-y)T^{-s/d})p_1((y-z)T^{-u/d})\\
&&\times\varphi(y)\varphi(z)\chi(T^{s-1})\chi(T^{s-1}+T^{u-1})\frac{1}{1+|x|^d}dudsdydzdx.
\end{eqnarray*}
In the integral $\int dx$ we pass to polar coordinates 
$(r,w) \;(r=|x|)$ and then substitute $r'=r^d$.

\noindent
We have
\begin{eqnarray}
\lefteqn{I_1(T)}\nonumber\\
&=&\frac{1}{\log T}\frac{1}{d}\int^T_1\int_{S_{d-1}}\int_{\errita^{2d}}\int^1_0\int^1_0\UNO_{[0,\frac{\log(T-1)}{\log T}]}(s)
\UNO_{[0,\frac{\log(T-T^s)}{\log T}]}(u)p_1((wr^{1/d}-y)T^{-s/d})p_1((y-z)T^{-u/d})\nonumber\\
\label{eq:3.58}
&&\times\varphi(y)\varphi(z)\chi(T^{s-1})\chi(T^{s-1}+T^{u-1})\frac{1}{1+r}dudsdydz\sigma_{d-1}(dw)dr,
\end{eqnarray}
where $\sigma_{d-1}$ is the Lebesgue measure on the unit sphere $S_{d-1}$ in $\erre^d$. Again, we use (\ref{eq:3.44}) putting $r'=\log r/\log T$, then it is easy to see that the integrand converges to
$$
\UNO_{[0,1]}(s)\UNO_{[0,1]}(u)\UNO_{[0,s]}(r)p^2_1(0)
\varphi(y)\varphi(z)\chi^2(0),
$$
and is bounded by $Cp^2_1(0)\varphi(y)\varphi(z)$. Hence (\ref{eq:3.54}) follows. By (\ref{eq:3.49}) and (\ref{eq:3.54})-(\ref{eq:3.57}) we obtain (\ref{eq:3.13}).

Now we pass to the proof of (\ref{eq:3.14}). By (\ref{eq:3.12}) and (\ref{eq:3.9}) for $\Phi$ of the form (\ref{eq:3.15}), after obvious substitutions we have
$$I\!I(T)\leq\frac{C}{(\log T)^{ 9/2}}\int_{\errita^d}\int^T_0{\cal T}_s\left(
\varphi\int^T_0{\cal T}_u\left(\varphi\int^T_0{\cal T}_r\varphi dr\right)du\right)(x)ds\frac{1}{1+|x|^d}dx.
$$

Using (\ref{eq:3.41}) twice we get
\begin{eqnarray*}
I\!I(T)&\leq&\frac{C}{(\log T)^{5/2}}
\int_{\errita^d}\int^T_0{\cal T}_s\varphi(x)ds\frac{1}{1+|x|^d}dx\\
&=&A+B,
\end{eqnarray*}
where
$$A=\frac{C}{(\log T)^{5/2}}\int_{|x|^d>T}\ldots\leq
\frac{C_1}{(\log T)^{3/2}}$$
by (\ref{eq:3.43}), and
\begin{eqnarray}
B&=&\frac{C}{(\log T)^{5/2}}\int_{|x|^d\leq T}\ldots\leq
\frac{C_2}{(\log T)^{3/2}}\int_{|x|^d\leq T}\frac{1}{1+|x|^d}dx,\\
&\leq&\frac{C_3}{(\log T)^{1/2}}.\nonumber
\end{eqnarray}
In the first estimate in (3.60) we have used (\ref{eq:3.41}) once more. Hence (\ref{eq:3.14}) follows.

Passing to the proof of tightness, first observe that the method employed in the proof of (\ref{eq:3.3}) can be also used to obtain convergence of finite dimensional distributions of $X_T$. This fact has been already used in \cite{BGT4}; here we repeat briefly the argument. For $\varphi_1,\varphi_2,\ldots,\varphi_k\in{\cal S}(\erre^d)$, all $\varphi_j>0$, and $0\leq t_1\leq\ldots\leq t_k\leq 1$,
it is easy to see that $E\,{\rm exp}\{-\sum^k_{j=1}\langle X_T(t_j),\varphi_j\rangle\}$ has the form (\ref{eq:3.6}) with
$$\Psi(x,t)=\sum^k_{j=1}\varphi_j(x)\UNO_{[0,t_j]}(t),$$
and the corresponding $v_T$ given by (\ref{eq:3.7}).

Approximating $\Psi$ by smooth functions we obtain that (\ref{eq:3.8}) holds, and then we argue as before.

In particular $X_T(\varepsilon)$ converges in law.
Therefore, to prove tightness of $X_T$ in $C([\varepsilon,1],{\cal S}'(\erre^d))$ it suffices to show that $\langle X_T-X_T(\varepsilon),\varphi\rangle\Rightarrow 0$ in $C([\varepsilon, 1])$, for any $\varphi\in {\cal S}(\erre^d), \varphi\geq 0$.

Denote
$$w_T(t)=\frac{1}{F_T}\int^{t T}_{\varepsilon T}\langle N_s,\varphi\rangle ds,\quad t\geq \varepsilon.$$
By (1.2) it is clear that it is enough to show that $w_T$ and $Ew_T$ converge to $0$ in law in $C([\varepsilon,1])$. Both processes are increasing, so it suffices to prove that
$$\lim_{T\to\infty}Ew_T(1)=0.$$
By (\ref{eq:3.0}) and  substitution $x'=xs^{-1/\alpha}$,
\begin{eqnarray*}
Ew_T(1)&=&\frac{1}{F_T}\int^T_{\varepsilon T}\int_{\errita^{2d}}
p_1(x-ys^{-1/\alpha})\varphi(y)\frac{1}{1+|x|^ds}dxdyds\\
&=&J_1(T)+J_2(T),
\end{eqnarray*}
where
\begin{eqnarray*}
J_1(T)&=&\frac{1}{F_T}\int^T_{\varepsilon T}\int_{|x|\leq 1}\int_{\errita^d}\ldots,\\
J_2(T)&=&\frac{1}{F_T}\int^T_{\varepsilon T}\int_{|x|>1}\int_{\errita^d}\ldots .
\end{eqnarray*}
We have
\begin{eqnarray*}
J_2(T)&\leq&\frac{C}{F_T}\int^T_{\varepsilon T}s^{-1}ds=C
\frac{\log(1/\varepsilon)}{F_T}\to 0,\\
J_1(T)&\leq&\frac{C_1}{F_T}\int^T_{\varepsilon T}
\int_{|x|\leq 1}\frac{1}{1+|x|^ds}dxds\\
&=&\frac{C_2}{F_T}\int^T_{\varepsilon T}\int^1_0
\frac{1}{1+rs}drds\leq
\frac{C_3\log (1/\varepsilon)}{(\log T)^{1/2}}\to 0.
\end{eqnarray*}$\hfill\Box$

\subsection{Systems with more general intensity measures $\mu$}

In this section we consider a measure  $\mu$ of the form (\ref{eq:2.13}). We sketch the proof of Proposition 2.9 and we discuss Example 2.10.
\vglue.5cm
\noindent
{\bf Proof of Proposition 2.9} We concentrate on the case $\gamma<d<\alpha$. The other cases will be mentioned later.

First  notice that it suffices to assume that $\nu\equiv 0$ in 
(\ref{eq:2.13}), since it is easy to see that with our normalization the terms corresponding to $\nu$ will vanish in the limit. We repeat the steps of the proof of Theorem 2.1. Observe that boundedness of $h$ implies that (\ref{eq:3.32}) also holds in the present case, hence (\ref{eq:3.14}) is obtained in the same way as before. Also, the tightness is proved without any change.

It remains to show (\ref{eq:3.13}). Instead of (\ref{eq:3.25}) we have
\begin{equation}
\label{eq:3.59}
I(T)=\int_{\errita^d}g_T(x)
\frac{T^{\gamma/\alpha}}
{1+|xT^{1/\alpha}|^\gamma} h
(T^{1/\alpha}x)dx,
\end{equation}
where $g_T$ is defined by (\ref{eq:3.26}). We write
\begin{equation}
\label{eq:3.60}
I(T)=I_1(T)+I_2(T),
\end{equation}
where
\begin{equation}
\label{eq:3.61}
I_1(T)=\int_{\errita^d}g_\infty(x)
\frac{1}{|x|^\gamma}h(T^{1/\alpha}x)dx,
\end{equation} 
with $g_\infty$ given by (\ref{eq:3.27}), and
$$I_2(T)=
\int_{\errita^d}\left(g_T(x)\frac{T^{\gamma/\alpha}}
{1+|xT^{1/\alpha}|^\gamma}-
\frac{g_\infty(x)}{|x|^\gamma}\right)h(T^{1/\alpha}x)dx.$$
(\ref{eq:3.31}) implies that $\lim\limits_{T\to\infty}I_2(T)=0$.

Note that from assumption (\ref{eq:2.14}) it follows that
\begin{equation}
\label{eq:3.62}
\lim_{R\to \infty}\int_{\errita^d}a(x)h(Rx)dx=C\int_{\errita^d}a(x)dx
\end{equation}
for $a(x)=\UNO_{|x|\leq r}$, where $C$ is the limit (\ref{eq:2.14}). Hence, it is easy to see that (\ref{eq:3.62}) also holds for any symmetric integrable function $a$. The function $g_\infty(x)|x|^{-\gamma}$ is obviously symmetric and integrable $(\gamma<d$ and (\ref{eq:3.31})). Therefore (\ref{eq:3.62}) implies 
(\ref{eq:3.13}). This completes the proof in the case $\gamma<d<\alpha$.

In the remaining cases for $\gamma<d$, tightness and (\ref{eq:3.14}) follow immediately from the cooresponding proofs for $\mu$ of the form (\ref{eq:2.1}). Also, to obtain (\ref{eq:3.13}) we repeat the same steps, obtaining $I(T)$ in an analogous form as in (\ref{eq:3.59}). And then we apply (\ref{eq:3.62}).
$\hfill\Box$
\vglue.5cm
For $d=\gamma$ the method described above cannot be applied. For example, in the case $d=\gamma=1<\alpha$ the function $g_\infty(x)|x|^{-d}$ is not integrable. To prove (\ref{eq:3.13}) we would need existence of the limit
$$\lim_{T\to\infty}\frac{1}{\log T}
\int_{|x|\leq 1}g_\infty(x)
\frac{T^{1/\alpha}}{1+|xT^{1/\alpha}|}
h(xT^{\frac{1}{\alpha}})dx.$$
This limit is easy to obtain for $h\equiv 1$, but it is not clear how to formulate an elegant condition assuring its existence in a more general case.

The case $\gamma=d=\alpha$ is even  more complicated because in (\ref{eq:3.58}) we would have $h(wr^{1/\alpha})$ under the integrals.
\vglue.5cm
\noindent
{\bf Proof of non-existence of the limit in Example 2.10} It is obvious that the only nontrivial normalization is that given by (\ref{eq:2.2}), since $h$ is bounded and separated from $0$. Analogously as in the previous proof, convergence of $X_T$ is equivalent to convergence of $I_1(T)$ defined by (\ref{eq:3.61}). We will show that $I_1(T)$ does not converge. Let 
$T_n=n^{n\alpha+{\alpha/2}}, n=2,3,\ldots$. On the set 
$\{x:\frac{1}{\sqrt{n}}\leq |x|\leq \sqrt{n}\}$ we have
$$h(T_n^{1/\alpha}x)=u(n)=\left\{\begin{array}{llll}
1&{\rm if}&n&\hbox{\rm is even},\\
2&{\rm if}&n&\hbox{\rm is odd}.
\end{array}\right.
$$
It is clear that $\lim\limits_{n\to\infty}(I_1(T_n)-I'_1(n))=0$, where
$$I'_1(n)=u(n)
\int_{\frac{1}{\sqrt{n}}\leq|x|\leq\sqrt{n}}g_\infty(x)|x|^{-\gamma}dx,$$
and obviously $I'_1(n)$ does not converge.

\subsection{The finite measure case}

\vglue.5cm
\noindent
{\bf Proof of Lemma 2.12} We improve slightly the proof of Lemma 2 of Bingham
 \cite{Bi}. Let
\begin{equation}
\label{eq:3.63}
\alpha '=1-\frac{1}{\alpha}.
\end{equation}
It is easy to see using self-similarity that
\begin{equation}
\label{eq:3.64}
\lim_{\theta \rightarrow 0}\theta^{\alpha '}\int^\infty_0\int_{\errita} 
e^{-\theta s}\varphi(y)p_s(x-y)dyds=K\int_{\errita}\varphi(x)dx
\end{equation}
for any $\varphi\in{\cal S}(\erre)$, but in general the convergence is not uniform in $x\in\erre$ if $\varphi$ is not compactly supported. Therefore condition $(A)$ of Darling-Kac \cite{DK} is not satisfied, so, unlike Bingham, we cannot apply directly their theorem.
We prove that
\begin{equation}
\label{eq:3.65}
Z_T{\flechita}KL\lambda  
\end{equation}
(${\flechita}$ denotes convergence of finite dimensional distibutions), where $L$ is a continuous increasing process whose inverse is an $\alpha'$-stable subordinator. On the other hand, 
it is known (see \cite{Be}, Prop. 4, Ch. V; see also \cite{FG0}) that such  $L$ is the local time process at $0$ of $\zeta$. 

It is clear that in order to prove (\ref{eq:3.65}) it suffices to show that
\begin{equation}
\label{eq:3.66}
(\langle Z(t_1),\varphi_1\rangle,\ldots,\langle Z(t_k),\varphi_k\rangle)\Rightarrow K\left(L(t_1)\int_{\errita}\varphi_1(x)dx,\ldots,L(t_k)\int_{\errita}\varphi_k(x)dx\right)
\end{equation}
for any $t_1,\ldots,t_k\in [0,1],\varphi_1,\ldots,\varphi_k\in{\cal S}(\erre),\varphi_1,\ldots,\varphi_k\geq 0, k=1,2,\ldots$.

Fix $\varphi_1,\ldots,\varphi_k$ as above $(\varphi_j\not\equiv 0)$. Let $M_T$
be the measure on $\erre^k_+$ such that
\begin{equation}
\label{eq:3.67}
M_T([0,t_1]\times\ldots\times[0,t_k])=\frac{1}{K^k}E\prod^k_{j=1}\frac{\langle Z_T(t_j),\varphi_j\rangle}{\langle\lambda,\varphi_j\rangle}.
\end{equation}
(\ref{eq:3.66}) will be proved if we show
\begin{equation}
\label{eq:3.68}
\lim_{T\to\infty}M_T([0,t_1]\times\ldots\times[0,t_k])=EL(t_1)\ldots L(t_k)
\end{equation}
for all $t_1,\ldots,t_k\in[0,1]$. To this end, by Lemma 3 of \cite{Bi} it suffices to prove that
\begin{eqnarray}
\lefteqn{\lim_{T\to\infty}\int_{\errita^k_+}e^{-(\theta_1t_1+\ldots+\theta_kt_k)}M_T(dt_1,\ldots,dt_k)}\nonumber\\
\label{eq:3.69}
&=&\sum_\pi\left[(\theta_{\pi(1)}+\ldots+\theta_{\pi(k)})(\theta_{\pi(2)}+\ldots+\theta_{\pi(k)})\ldots\theta_{\pi(k)}\right]^{-\alpha'}
\end{eqnarray}
for all $\theta_1,\ldots,\theta_k\geq 0,$ the sumation being over all permutations $\pi$ of $\{1,\ldots,k\}$. For simplicity we will show (\ref{eq:3.69}) for $k=2$. Without loss of generality we may assume that
$\langle\lambda,\varphi_j\rangle =1, j=1,2$.

We have
\begin{equation}
\label{eq:3.70}
\int^\infty_0\int^\infty_0e^{-\theta t_1-\theta t_2}M_T(dt_1,dt_2)=\frac{1}{K^2}\left(J_1(T)+J_2(T)\right),
\end{equation}
where
\begin{eqnarray}
\label{eq:3.71}  
J_1(T)&=&\frac{1}{F^2_T}\int^\infty_0\int^{t_2}_0e^{-\theta_1t_1-\theta_2t_2}T^2E\varphi_1(\zeta_{Tt_1})\varphi_2(\zeta_{Tt_2})dt_1dt_2,\\
\label{eq:3.72}
J_2(T)&=&\frac{1}{F^2_T}\int^\infty_0\int^{t_1}_0\ldots dt_2dt_1,
\end{eqnarray}
where $\ldots$ denotes the same integrand as in (\ref{eq:3.71}).

By the Markov property, self-similarity and (\ref{eq:2.17}), we have
\begin{eqnarray*}
J_1(T)&=&\int^\infty_0\int^{t_2}_0\int_{\errita}\int_{\errita} 
e^{-\theta_1t_1-\theta_2t_2}\varphi_1(y)\varphi_2(z)t_1^{-1/\alpha}
(t_2-t_1)^{-1/\alpha}\\
&&\times p_1(T^{-1/\alpha}t_1^{-1/\alpha}y)p_1(T^{-1/\alpha}
(t_2-t_1)^{-1/\alpha}(z-y))dzdydt_1dt_2,
\end{eqnarray*}
hence
\begin{eqnarray*}
\lim_{T\to\infty}J_1(T)&=&\int^\infty_0\int^{t_2}_0
\frac{e^{-\theta_1t_1}}{t_1^{1/\alpha}}
\frac{e^{-\theta_2t_2}}{(t_2-t_1)^{1/\alpha}}
dt_2dt_1p^2_1(0)\\
&=&K^2(\theta_1+\theta_2)^{-\alpha'}\theta_2^{-\alpha'},
\end{eqnarray*}
where $K=\frac{1}{\pi\alpha}\Gamma(\frac{1}{\alpha})\Gamma(1-\frac{1}{\alpha}).$ The limit of $J_2(T)$ is calculated identically, so we obtain 
(\ref{eq:3.69}).

For $\varphi\geq 0$, (\ref{eq:3.65}) implies that
$\langle Z_T,\varphi\rangle\Rightarrow KL\langle\lambda,\varphi\rangle$ in $C([0,1])$, since $\langle Z_T,\varphi\rangle$ is an increasing process. 
From this it follows immediately that $\{\langle Z_T,\varphi\rangle\}_T$ is tight for any $\varphi\in{\cal S}(\erre^d$). Hence the proof of the lemma is complete by Mitoma's theorem \cite{M}. $\hfill\Box$
\vglue.5cm
\noindent
{\bf Proof of Theorem 2.11} It is easy to see that Lemma 2.12 remains true with the same limit if $\zeta$ is replaced by $x+\zeta$.

On the other hand,
$$\langle Y_T(t),\varphi\rangle=\sum_{x_j\in N_0}\frac{1}{F_T}\int^{tT}_0\varphi(x_j+\zeta^j_s)ds,$$
where $\zeta^1,\zeta^2,\ldots$ are independent copies of $\zeta$, independent of $N_0$. Now the theorem follows from Lemma 2.12 and the fact that $N_0(\erre)$ has the same law as $\nu$. $\hfill\Box$

\vglue.5cm
\noindent
{\bf Proof of Theorem 2.13} Let $Z_T$ be defined by (\ref{eq:2.19}) with $F_T=\log T$. It suffices to prove that
\begin{equation}
\label{eq:3.73}
\langle Z_T(1),\varphi\rangle\Rightarrow p_1(0)
\rho_1\langle\lambda,\varphi\rangle
\end{equation}
for $\varphi\geq 0$. Indeed, (\ref{eq:3.73}) implies that $\frac{1}{\log T}\int^{t_2T}_{t_1T}\varphi(\zeta_s)ds$ converges to $0$ in probability for any $0<t_1\leq t_2$. Now it is easy to see that
$$(\langle Z_T(t_1),\varphi_1\rangle,\ldots,\langle Z_T(t_k),\varphi_k\rangle)\Rightarrow p_1(0)(\rho_1\langle\lambda,\varphi_1\rangle,\ldots,
\rho_1\langle\lambda,\varphi_k\rangle)$$
for $0<t_1\leq\ldots\leq t_k$. Hence, analogously as before we obtain that 
$Z_T\Rightarrow K\rho_1\lambda$ in $C([\varepsilon,1],{\cal S}'(\erre^d))$, and this easily implies (\ref{eq:2.20}).

Observe that (\ref{eq:3.73}) has exactly the form as in the Darling-Kac theorem
\cite{DK}, but we cannot apply it directly, since
$$\frac{1}{-\log\theta}\int^\infty_0\int_{\errita^d} e^{-\theta s}\varphi(y)p_s(x-y)dyds$$
may not converge uniformly in $x$ as $\theta\to 0$ (so,  condition (A) is not satisfied). Nevertheless, the proof of the Darling-Kac theorem can be repeated with some care in the present case yielding the desired result. 
$\hfill\Box$
\vglue.5cm
\noindent
{\bf Appendix}
\noindent
\vglue.5cm
\noindent
{\bf Proposition} {\it For $\gamma>0$,}
$$E\langle N_t,\varphi\rangle\to 0\quad{ as}\quad 
t\to\infty, \quad\varphi\in {\cal S}(\erre^d).$$
\vglue.5cm
\noindent
{\bf Proof} It suffices to assume $\varphi\geq 0.$

Let $\gamma>d$. Using self-similarity of the standard $\alpha$-stable process, we have from (3.1),
\begin{eqnarray*}
\lefteqn{E\langle N_t,\varphi\rangle=\int_{\errita^d}
\int_{\errita^d}p_t(x-y)\varphi(y)dy\frac{dx}{1+|x|^\gamma}}\\
&&=t^{-d/\alpha}\int_{\errita^d}\int_{\errita^d}p_1(t^{-1/\alpha}(x-y))\varphi(y)dy\frac{dx}{1+|x|^\gamma},
\end{eqnarray*}
and the integrals converge to a finite limit
as $t\to\infty$.

Let $\gamma\leq d$. By the following version of 
Young's inequality \cite{LL} (Theorem 4.2),
$$\left|\int_{\errita^d}\int_{\errita^d}f(x)g(x-y)h(y)dxdy\right|\leq C||f||_p||g||_q||h||_r,$$%
$ p,q,r\geq 1, 1/p+1/q+1/r=2$, we have
$$
E\langle N_t,\varphi\rangle
\leq C\left(
\int_{\errita^d}p_t(x)^pdx
\right)^{1/p}
\int_{\errita^d}\varphi(x)dx
\left(\int_{\errita^d}\frac{1}{(1+|x|^\gamma)^r}dx
\right)^{1/r},
$$
where $p>1$, $q=1, r>{d}/{\gamma}\geq 1$. The integrals are finite. Again by self-similarity, and the change of variable $x=t^{1/\alpha}y$,
$$\left(\int_{\errita^d}p_t(x)^pdx\right)^{1/p}=t^{-(1-1/p)d/\alpha}
\left(\int_{\errita^d}p_1(y)^pdy\right)^{1/p},$$
where $1-1/p>0$. 

In both cases the result is obtained. 

$\hfill\Box$
\vglue.5cm
\noindent
{\bf Acknowledgment}
\vglue.5cm
We thank the hospitality of the Institute of Mathematics, National University of Mexico (UNAM), where this paper was partially written.

%\normalfont
\renewcommand{\refname}


{\normalsize References}
\begin{thebibliography}{99999}
%\renewcommand{\refname}{\small References}

\bibitem[1]{Be}  J. Bertoin, L\'evy Processes, Cambridge Univ. Press, 1996.
\bibitem[2]{Bi}  N.H. Bingham, Limit theorems for occupation times of Markov processes, Z. Wahrschein. verw. Geb. 17 (1971) 1-22.
\bibitem[3]{BGR} T. Bojdecki, L.G. Gorostiza and S. Ramaswami, Convergence of ${\cal S}'$-valued processes and space-time random fields, J. Funct. Anal. 66 (1986) 21-41.
\bibitem[4]{BGT0}T. Bojdecki, L.G. Gorostiza, A. Talarczyk, Fractional Brownian density process and its self-intersection local time of order $k$, J. Theor. Probab. 17 (2004) 717-739.
\bibitem[5]{BGT1}T. Bojdecki, L.G. Gorostiza, A. Talarczyk, Limit theorems for occupation time fluctuations of branching systems I: Long-range dependence, Stoch. Proc. Appl. 116 (2006) 1-18. 
\bibitem[6]{BGT2} T. Bojdecki, L.G. Gorostiza, A. Talarczyk, Limit theorems for occupation time fluctuations of branching systems II: Critical and large dimensions, Stoch. Proc. Appl. 116 (2006) 19-35.
\bibitem[7]{BGT3} T. Bojdecki, L.G. Gorostiza, A. Talarczyk, A long range dependence stable process and an infinite variance branching system, Ann. Probab.
 (to appear). Math. ArXiv PR/0511739.
\bibitem[8]{BGT4}T. Bojdecki, L.G. Gorostiza, A. Talarczyk, Occupation time fluctuations of an infinite variance of branching system in large dimensions, Bernoulli (to appear). Math. ArXiv PR/0511745.
\bibitem[9]{BGT5} T. Bojdecki, L.G. Gorostiza, A. Talarczyk, Some extensions of fractional Brownian motion and sub-fractional Brownian motion related to particle systems (in preparation).
\bibitem[10]{CG} J.T. Cox, D. Griffeath, Large deviations for Poisson systems of independent random walks, Z. Wahrschein. verw. Geb. 66 (1984) 543-558.
\bibitem[11]{DK} D.A. Darling, M. Kac, On occupation times for Markoff processes, Trans. Amer. Math. Soc. 84 (1957) 444-458.
\bibitem[12]{DGW} D.A. Dawson, L.G. Gorostiza, A. Wakolbinger, Occupation time fluctuations in branching systems, J. Theor. Probab. 14 (2001) 729-796.
\bibitem[13]{DR}  J.-D. Deuschel, J. Rosen, Occupation time large deviations for critical branching Brownian motion, super-Brownian motion and related processes, Ann. Probab. 26 (1998) 602-643.
\bibitem[14]{DW}  J.D. Deuschel, K. Wang, Large deviations for the occupation time functional of a Poisson system of independent particles, Stoch. Proc. Appl. 52 (1994) 183-209.
\bibitem[15]{FG0}P.J. Fitzsimmons, R.K. Getoor, On the distribution of the Hilbert transform of the local time of a symmetric L\'evy process, Ann. Probab. 20 (1992) 1484-1497.
\bibitem[16]{FG1} P.J. Fitzsimmons, R.K. Getoor, Limit theorems and variation properties for fractional derivatives of the local time of a stable process, Ann. Inst. H. Poincar\'e, Probab. Stat. 28 (1992) 311-333.
\bibitem[17]{GW}  L.G. Gorostiza, A. Wakolbinger, Persistence criteria for a class of critical branching particle systems in continuous time, Ann. Probab. 19 (1991) 266-288.
\bibitem[18]{HS}  R.A. Holley, D.W. Stroock, Generalized Ornstein-Uhlenbeck processes and infinite particle branching Brownian motions, Publ. Res. Inst. Math. Sci. 14 (1978) 741-788.
\bibitem[19]{I} I. Iscoe, A weighted occupation time for a class of measure-valued branching processes, Probab. Theor. Relat. Fields 71 (1986) 85-116.
\bibitem[20]{K} A. Klenke, Multiple scale analysis of clusters in spatial branching models, Ann. Probab. 25 (1997) 1670-1711.
\bibitem[21]{LL}E.H. Lieb, M. Loss, Analysis, 2nd. Edition, Amer. Math. Soc., Providence, 2001.
\bibitem[22]{ML} A. Martin-L\"of, Limit theorems for the motion of a Poisson system of independent Markovian particles with high density, Z. Wahrschein. 
verw. Geb. 34 (1976) 205-223.
\bibitem[23]{Mi} P. Mi\l o\'s, Occupation time fluctuations of Poisson and equilibrium finite variance branching systems, Prob. Math. Stat. (to appear). Math. ArXiv PR/0512414.
\bibitem[24]{M} I. Mitoma, Tightness of probabilities in $C([0,1], {\cal S}')$ and $D([0,1], {\cal S}')$, Ann. Probab. 11 (1983), 989-999. 
\bibitem[25]{T} A. Talarczyk, A functional ergodic theorem for the occupation time process of a branching system (preprint).
\bibitem[26]{W} J.B. Walsh, An introduction to stochastic partial differential equations, Ecole d'Et\'e de Probabilit\'es de Saint-Flour XIV-1984. Lect. Notes Math. 1180, Springer, Berlin, 265-439.
\end{thebibliography}
\end{document}